\documentclass[graybox]{svmult}

\usepackage{type1cm}        
\usepackage{makeidx}         
\usepackage{graphicx}        
\usepackage{multicol}        
\usepackage[bottom]{footmisc}

\usepackage{newtxtext}       %
\usepackage[varvw]{newtxmath}       

\makeindex             

\begin{document}

\title*{Efficient implicit solvers for models of neuronal networks}
\author{Luca Bonaventura\orcidID{0000-0002-1994-0217},\\ Soledad Fernández-García\orcidID{0000-0001-6993-407X} and \\ Macarena Gómez-Mármol\orcidID{0000-0003-2651-5689} \\}
\institute{Luca Bonaventura \at Politecnico di Milano, Via Bonardi 9, 20133 Milano, Italy, \email{luca.bonaventura@polimi.it}
\and Soledad Fernández-García \at Universidad de Sevilla, Apdo.\,de correos 1160,  41080 Sevilla, Spain \email{soledad@us.es}
\and Macarena Gómez-Mármol \at Universidad de Sevilla, Apdo.\,de correos 1160,  41080 Sevilla, Spain \email{macarena@us.es}}
%
%
\maketitle

\abstract*{We introduce economical versions of standard implicit ODE solvers that are specifically tailored for the efficient and accurate
simulation of neural networks. The specific versions of the ODE solvers proposed here, allow to achieve a significant increase in the efficiency
 of network simulations, by reducing the size of the algebraic system being solved at each time step, a technique inspired by  
 very successful semi-implicit approaches in computational
fluid dynamics and structural mechanics.
While we focus here specifically on Explicit first step, Diagonally Implicit Runge Kutta methods (ESDIRK), 
 similar simplifications can also be applied to any implicit ODE solver. In order to demonstrate the capabilities of the proposed methods,
we consider networks based on three different single-cell models with   slow-fast dynamics,
including  the classical FitzHugh-Nagumo model, a Intracellular Calcium Concentration model and the Hindmarsh-Rose model.
Numerical experiments on the simulation of networks  of increasing size based on these models demonstrate the increased efficiency of the proposed methods.}

\abstract{We introduce economical versions of standard implicit ODE solvers that are specifically tailored for the efficient and accurate
simulation of neural networks. These reformulations allow to achieve a significant increase in the efficiency
 of network simulations, by reducing the size of the algebraic systems effectively solved at each time step.
While we focus here specifically on Explicit first step, Diagonally Implicit Runge Kutta methods (ESDIRK), 
 similar simplifications can also be applied to any implicit ODE solver. In order to demonstrate the capabilities of the proposed methods,
we consider networks based on three different single-cell models with   slow-fast dynamics,
including  the classical FitzHugh-Nagumo model, a Intracellular Calcium Concentration model and the Hindmarsh-Rose model.
Numerical experiments on the simulation of networks  of increasing size based on these models demonstrate the  superior efficiency of the proposed
economical methods.}

\section{Introduction}
 \label{sec:intro} \indent
 
Synchronization between neuronal activities plays an important role in the understanding of the nervous system. Starting with the seminal work of Hodgkin-Huxley \cite{hodgkin:1952}, the complexity of the ionic dynamics is usually reflected in the models of
neural activity by the presence of nonlinearities and of different timescales for the different variables. 
The rich variety of synchronization types
that can take place in neuron networks results both from the complexity of the neural  dynamics and
from the    scale and structure of the network itself, which can vary  from a small number of cells
(microscopic scale), through neuron populations (mesoscopic scale), to large areas
of the brain and spinal cord (macroscopic scale).

The synchronization properties of coupled systems with multiple timescales, such as relaxation oscillators \cite{fitzhugh:1961, nagumo:1962}, bursters \cite{hindmarsh:1984}  and systems presenting Mixed-Mode Oscillations (MMOs)  \cite{desroches:2012}, differ strongly from those
between coupled harmonic oscillators, and the role of the coupling strength is
different in the two cases. In particular, it has been shown  that canard phenomena \cite{benoit:1981}
arising in multiple timescale systems play a prominent role in organizing the
synchronization of coupled slow-fast systems. 
Recent studies  on these topic mostly address
issues such as synchronization and desynchronization,
local oscillations and clustering \cite{ermentrout:2001} and the
problem of synchronization of coupled multiple timescale systems  constitutes a very active field of research. Furthermore,
neuronal networks  with similar properties are the main component of many supervised learning methods based on recurrent neural networks,
such as continuous time Liquid State Machines, see e.g.,  \cite{maass:2002, maass:2011}
and Echo State Networks \cite{jaeger:2007, yildiz:2012}.

Due to the nonlinearities involved and to the scale of the networks under study, numerical simulation is an essential tool to understand
and simulate synchronization phenomena.
From a numerical point of view, efficient simulations of neural networks require the use of special methods suitable for stiff problems, due to the slow-fast nature of the dynamics. Furthermore, if the number of cells in the cluster is large, numerical simulations can entail a substantial computational effort if  standard ODE solvers are applied. Numerical techniques with similar properties are also required in the so-called Neural ODE approaches to neural network modelling \cite{chen:2018}.

In this work,  we  show how to build economical versions of standard implicit ODE solvers specifically tailored for the efficient and accurate
simulation of neural networks. The specific versions of the ODE solvers proposed here  allow to achieve a significant increase in the efficiency
 of network simulations, by reducing the size of the algebraic system being solved at each time step.  This development is inspired by  
 very successful semi-implicit approaches in computational
fluid dynamics, see e.g. \cite{bonaventura:2000}. A similar approach was applied
in \cite{bonaventura:2021a} to the classical equations of structural mechanics.
While we focus here specifically on Explicit first step, Diagonally Implicit Runge Kutta methods (ESDIRK), 
see the reviews \cite{kennedy:2016, kennedy:2019}, analogous simplifications can be applied to any
implicit ODE solver. 

In order to demonstrate the capabilities of the proposed methods,
we consider networks based on three different single-cells models, aiming to cover a wide variety of models, with different dynamical properties depending, among others, on their slow-fast nature. Notice that the slow and fast label only concern the different time scales arising in the equations. Due to the nonlinear coupling present in all the considered models, both kinds of variables experience sharp transients which need to be accurately simulated by numerical methods. The first model is the classical FitzHugh-Nagumo (FN) system \cite{fitzhugh:1961,nagumo:1962}. It consists on a system of two equations that evolve with different time scales. With only two equations, the system is able to reproduce the neuron excitability. 
The second model is a FN system with an extra variable representing the Intracellular Calcium Concentration (ICC) in neurons \cite{krupa:2013}. 
The third variable is slow, so that the resulting system is a two slow - one fast system. This allows the system to have MMOs, that is, oscillatory patterns with an alternation of small and large amplitude oscillations \cite{desroches:2012}.  
The third model is the Hindmarsh-Rose (HR) system \cite{hindmarsh:1984}. This is a system with one slow and two fast variables, 
which displays bursting oscillations.  The main characteristic of these oscillations is an
alternation of slow phases, where the system is
quasi-stationary, and rapid phases, where the system is
quasi-periodic. During the latter phase, the system solutions display groups of large-amplitude
 oscillations or spikes that occur on a faster timescale  \cite{izhikevich:2000, rinzel:1987}.
 For classical parameter values, the HR model produces
square-wave bursting, one of the three main classes introduced in \cite{rinzel:1987}, taking this name because of the form of the oscillations.  

From each of these   single-cell models, following the approaches proposed in the literature \cite{campbell:2001,eteme:2017,fernandez:2020,ibrahim:2019,krupa:2013,yong:2008}, we build three different networks, 
based on the density of the  coupling matrix (sparse, middle, dense), with the objective of testing the efficiency of the developed methods 
in these three different situations.
 
The rest of the article is outlined as follows: in Section \ref{sec:model}, we present the single neuron models that we use at each node of the networks. In Section \ref{nt:model}, we build the networks that we aim to simulate. Section \ref{sec:solvers} is devoted to  the derivation of efficient implicit solvers for the Implicit Euler method adapted to each network. After that, in Section \ref{sec:esdirk} we  discuss the extension of the procedure
 outlined in  Section \ref{sec:solvers}   to a class of convenient high order ODE solvers. Numerical experiments are performed in Section \ref{sec:numerical}. Finally, Section \ref{sec:conclu} is devoted to exposing conclusions and perspectives of the present work.

 \section{Single neuron models}
 \label{sec:model} \indent 
 
 We present here the single-cell neuron models that we use to construct the networks in Section \ref{nt:model}. As we have already commented in the Introduction, we consider three different slow-fast models: 
the classical FN system \cite{fitzhugh:1961,nagumo:1962}, the FN system with an extra variable representing the ICC in neurons \cite{krupa:2013}
and the  HR system \cite{hindmarsh:1984}. We consider first the classical FN model of spike generation \cite{fitzhugh:1961, nagumo:1962}, given by
 
\begin{equation*}\label{FNsystem}
\begin{array}{rcl}
      \dot x&=& -y+f(x),\\
      \dot y&=&\varepsilon (x+g(y)), \\
      \end{array}
\end{equation*}
where
\begin{equation}\label{def:f}
f(w)=4w-w^3, \  \   \   g(w)=a_1 w+a_2,
\end{equation}
with $(x,y)\in\mathbb{R}^2.$  Here $x$ represents the membrane potential and $y$ is the slow recovery variable. The timescale separation parameter $\varepsilon$ fulfills $0<\varepsilon\ll 1,$ and $a_1,a_2\in\mathbb{R}$ are usually taken such that the system has only one equilibrium point.

 We consider then the ICC model, introduced in \cite{krupa:2013},  which is given by
 \begin{equation}\label{ICCsystem}
\begin{array}{rcl}
      \dot x&=&\tau( -y+f(x)-\phi_{f}(z)),\\
      \dot y&=&\tau \varepsilon k(x+g(y)), \\
       \dot z&=& \tau\varepsilon \left( \phi_{r}(x)+r(z)\right),
      \end{array}
\end{equation}
where functions $f$ and $g$ are given in \eqref{def:f},

\begin{equation}\label{def_phi}
\phi_f(w)=\frac{\mu w}{w+z_0}, \quad \phi_r(w)=\frac{\lambda}{1+\exp{(-\rho(w-x_{on})})}, 
\end{equation}
\begin{equation}\label{def_r}
 r(w)=-\frac{w-z_b}{\tau_z},
 \end{equation}
and $(x,y,z)\in\mathbb{R}^3$.   Here, $x$ represents the membrane potential, $y$ is the slow recovery variable and $z$ stands for the ICC. The timescale separation parameter $\varepsilon$ fulfills $0<\varepsilon\ll 1. $ The
parameter $\tau>0$ has been introduced in \cite{krupa:2013} so that the outputs  comply with  a given physical timescale and does not impact the phase portrait. Moreover, following \cite{fernandez:2020}, we assume $a_1<0,$ $|a_1|\ll 1$ and parameters $a_2,z_0,\lambda,\tau_z,z_b,k$ to be strictly positive. We also consider a large enough $\rho$ value, so that the sigmoid $\phi_r$ is steep at its inflection point and represents a sharp activation function.  

Finally, we consider the HR model proposed in \cite{hindmarsh:1984} 
 \begin{equation*}\label{HRsystem}
\begin{array}{rcl}
      \dot x&=&l(x)+y-z+I,\\
      \dot y&=&c+m(x)-y, \\
       \dot z&=&\varepsilon(k(x-x_0)-z),
      \end{array}
\end{equation*}
where 
\begin{equation}\label{deflm}
l(w)=-aw^3+bw^2,\quad m(w)=-dw^2,
\end{equation}
 with $(x,y,z)\in\mathbb{R}^3$.  Here $x$ represents the membrane potential and  $y$ and $z$ take into account the transport of ions across the membrane through the ion channels. The timescale separation parameter $\varepsilon$ fulfills $0<\varepsilon\ll 1.$
Parameter $k>0$ has been introduced so that the outputs fit realistic biological evolution patterns and parameter $I$ stands for the external exciting current. Following \cite{eteme:2017}, we take  $I=3.28 $ and $k=4, $ $\varepsilon=0.008. $  Finally, $a,b,c$ and $d$ are positive parameters, with classical values  $a = 1, b = 3, c = 1,$ and $d = 5,$ and $x_0<0,$ usually taken as $x_0=-1.6.$

 \section{Network models}
 \label{nt:model} \indent
 
 In this Section  we  show how to construct the networks that will be simulated in this work, starting from the single-cell neuron models considered in Section \ref{sec:model}. 
The special structure of the Jacobians for the resulting models will then be exploited in Section \ref{sec:solvers} to derive implicit numerical methods that reduce substantially the computational cost by requiring the solution of linear systems of order $N$ rather than $2N$ or $3N $ in the Newton method iterations. Ultimately, the same ODE system will be solved, but a significant reduction in the number of required operations will be achieved by simplifying the nonlinear iterations required by implicit methods.

From now on, we denote by $\mathbb{I}_N$ the identity matrix of order $N,$ $\mathbb{O}_N$ the null matrix of order $N $ and $e$ the vector 
  $e=[1,1,...,1]^T \in \mathbb{R}^N.$ 
   
 \subsection{FitzHugh-Nagumo neuron network}
 \label{ssec:FHN_nn}
Coupled FN systems  have been widely considered in the literature, see, for instance,  \cite{campbell:2001, ibrahim:2019, yong:2008}.  As a first case, we consider the model given by model variables $ x, y,\in \mathbb{R}^N $ and defined,
with a slight abuse of notation that will be repeated for all the models considered in this paper,
 by the functions 
 
\begin{equation}\label{def:fyg}
 {f}( {w})=[f(w_1),\dots,f(w_N)]^T ,\quad  {g}( {w})=[g(w_1),\dots,g(w_N)]^T,
 \end{equation}
 where functions $f$ and $g$ are given in expressions \eqref{def:f}. Different cells are connected  through normalized electrical coupling, represented  by the 
 symmetric connectivity matrix  $\tilde{{C}} =(\tilde c_{ij}), $  with 
 $\tilde c_{ij}\in[-1,1], $  so that an undirected network is obtained.  The model equations are then  
   
  \begin{equation*}\label{FNnetwork}
   \begin{array}{l}
      \dot x_i= -y_i+f(x_i)+\dfrac{1}{N}\sum\limits_{j=1}^N\tilde c_{ij}(x_i-x_j), \\
      \dot y_i=\varepsilon\left(x_i + g(y_i)\right).\\
    \end{array}\
    \end{equation*}
  for $  i=1,2,...N. $ Setting also $c=(c_i)$ and $D=(d_{ij}),$ with
 $$ c_i=\frac{1}{N} \sum\limits_{j=1}^N\tilde c_{ij},  \ \ \ 
 d_{ij}=   c_i\delta_{ij} -\frac1N\tilde c_{ij}, $$ where $\delta_{ij} $ denotes the Kronecker delta,  
 one can then rewrite the equation for $x_i$ as  
 $$
  \dot x_i= f(x_i)-y_i +\sum\limits_{j=1}^N d_{ij}x_j.
  $$
  The whole system 
 can then be written in vector notation as

  \begin{eqnarray}
 \label{sys1}
 \left [\begin{array}{c}
\dot { {x}}   \\ 
 \dot { {y}}   
\end{array} \right ] &=&  \left [\begin{array}{c }
 {f}( {x}) - {y} +  {D} {x}   \\ 
 \varepsilon  {x}    +  \varepsilon  {g}( {y})  
\end{array} \right ] 
= \left [\begin{array}{cc}
 {D}& -{\mathbb{I}_N} \\
 \varepsilon{\mathbb{I}_N}   &   \varepsilon   
 a_1{\mathbb{I}_N}
 \end{array} \right ]\left[\begin{array}{c}
 {x}   \\ 
 {y}   \\  
\end{array} \right ] + \left [\begin{array}{c}
 {f}( {x})  \\ 
  \varepsilon a_2{e} \\ 
\end{array} \right ].
  \end{eqnarray}
 As a consequence, the Jacobian of the vector field on the right hand side is given by
 \begin{equation} 
 \label{j_sys1}
  {J}=  \left [\begin{array}{cc}
{D}+{\rm diag}({f}^{\prime}({x})) & -{\mathbb{I}_N}  \\
 \varepsilon  {\mathbb{I}_N} &  \varepsilon   a_1{\mathbb{I}_N}
 \end{array} \right].
 \end{equation}

\subsection{Intracellular Calcium Concentration neuron network}
 \label{ssec:ICC_nn}
 
We now build a neuron network by coupling systems of the form \eqref{ICCsystem}, introduced in \cite{krupa:2013}. A first step in the construction of a network from the one-cell system \eqref{ICCsystem} has been given in \cite{fernandez:2020}, where the synchronization patterns of the symmetric coupling of two cells were analyzed. 
Here, we  consider instead a larger and in principle arbitrary  network.  
 We consider the model  variables $ {x}, {y}, {z},\in \mathbb{R}^N $ and 
 model coefficients
 $\tau, \varepsilon,  \lambda, \mu, \tau_z, x_{on}, z_0,$
 defined by the functions $ {f}$ and $ {g}$ given in expression \eqref{def:fyg} and
 $\phi_f( {w})=[\phi_f(w_1),\dots,\phi_f(w_N)]^T,  $ 
 $ {\phi}_r( {w})=[\phi_r(w_1),\dots,\phi_r(w_N)]^T,$
$  {r}( {w})=[r(w_1),\dots,r(w_N)]^T, $
 where functions, $\phi_f$ and $\phi_r$ are given in \eqref{def_phi} and $r$ is given in \eqref{def_r}.
We then define  $ {k}=[k_1,\dots,k_N]^T $ as  random numbers uniformly distributed in the interval $[0.6,1.4].$ 
Different cells are connected by the  symmetric connectivity matrix  $\tilde{ {C}} =(\tilde c_{ij}), $  with 
 $\tilde c_{ij}\in[-1,1], $  so that an undirected network is obtained.  The model equations are then   

 \begin{equation*}\label{ICCnetwork}
   \begin{array}{l}
      \dot x_i= \tau(-y_i+f(x_i)-\phi_f(z_i)), \\
      \dot y_i=\tau\varepsilon k_i\left(x_i + g(y_i)+\frac{2}{N}\sum\limits_{j=1}^N\tilde c_{ij}(x_i-x_j)\right),\\
      \dot z_i= \tau \varepsilon \left(\phi_r(x_i)+r(z_i)\right).   \\
    \end{array}
    \end{equation*}
     for
 $i=1,2,...N. $ Note that, in contrast to the first case, we place the coupling in the second equation, following the phenomenological approach considered previously in \cite{fernandez:2020}. Furthermore, we multiply by $2/N$ because we consider two different clusters of neurons, as done also in \cite{bandera:2022} inspired by applications to motoneurons \cite{fallani:2015}. Setting also  $c=(c_i),$  $D=(d_{ij})$ and  $K=(k_{ij}),$ with
 $ c_i= \frac{2}{N} \sum\limits_{j=1}^N\tilde c_{ij}, $ 
  $d_{ij}=   (1+c_i)\delta_{ij} -2\frac{\tilde c_{ij}}{N}, $ 
 $  k_{ij}=\delta_{ij} k_i,$
 where $\delta_{ij} $ denotes the Kronecker delta,  the equation for $y_i$ becomes
 $$
  \dot y_i=\tau\varepsilon k_i\left( g(y_i)+     \sum\limits_{j=1}^Nd_{ij}x_j\right),
  $$
  and the whole system in vector notation as
  
 \begin{eqnarray}
 \label{sys2}
 \left [\begin{array}{c}
\dot{ {x}}   \\ 
 \dot{ {y}}   \\ 
 \dot{ {z}}   \\ 
\end{array} \right ]=\tau \left [\begin{array}{c }
 {f}( {x}) - {y} -  {\phi}_f( {z})  \\ 
 \varepsilon {K} {D} {x}    +  \varepsilon {K}  {g}( {y})  \\ 
\varepsilon {\phi}_r( {x})  +\varepsilon {r}( {z})     \\ 
\end{array} \right ] 
&=& \tau\left [\begin{array}{ccc}
{\mathbb{O}_N}
  & - {\mathbb{I}_N}& {\mathbb{O}_N}\\
 \varepsilon{K}{D}  & \varepsilon a_1\ {K} &   {\mathbb{O}_N}\\
   {\mathbb{O}_N}&   {\mathbb{O}_N}& -\varepsilon {\mathbb{I}_N}/\tau_z \\
 \end{array} \right ]\left[\begin{array}{c}
 {x}   \\ 
  {y}   \\ 
  {z}   \\ 
\end{array} \right ] 
\nonumber\\
&+&\tau\left [\begin{array}{c}
 {f}( {x}) -  {\phi}_f( {z}) \\ 
\varepsilon a_2 {k}   \\ 
\varepsilon {\phi}_r( {x}) + \varepsilon z_b/\tau_z    \\ 
\end{array} \right ].
 \end{eqnarray}
 As a consequence, the Jacobian of the vector field on the right hand side is given by
 \begin{equation} 
 \label{j_sys2}
  {J}= \tau  \left [\begin{array}{ccc}
  {\rm diag}({f}^{\prime}({x}))& - {\mathbb{I}_N}&-{\rm diag}({\phi}^{\prime}_f({z}))\\
 \varepsilon{K}{D}  & \varepsilon a_1{K} &   {\mathbb{O}_N}\\
 \varepsilon {\rm diag}({\phi}^{\prime}_r({x}))&   {\mathbb{O}_N}& -\varepsilon {\mathbb{I}_N}/\tau_z \\
 \end{array} \right ],
 \end{equation}
 where
 $$   \phi^{\prime}_f(w)=\frac{\mu z_0}{(w+z_0)^2}, \ \ \
  \phi^{\prime}_r(w)=\frac{\lambda \rho\exp{(-\rho(w-x_{on}))}}{ [ 1+\exp{(-\rho(w-x_{on}))} ]^2} .$$
 
   \subsection{Hindmarsh-Rose neuron network} 
   \label{ssec:HR_nn}
In this Section we consider instead the HR neuron network given in \cite{eteme:2017}, which consist on a network of electrically coupled HR systems.  
The model variables are ${x},{y},{z}\in \mathbb{R}^N .$
Different cells are connected through normalized electrical coupling  by a symmetric connectivity matrix  $  {{C}} =( c_{ij}).$ 
In \cite{eteme:2017}, a full connectivity matrix was considered given by
 $   c_{ij}=|i-j|^{-2} $ for $i\neq j=1,2,...N $ and $   c_{ij}=0 $ for $i=j=1,2,...N. $  
 However, alternative structures could also be considered, in particular given by sparse connectivity matrices.
 The model equations are then for
 $i=1,2,...N,$ 
 
 \begin{equation*}\label{HRnetwork}
\begin{array}{l}
      \dot x_i=l(x_i) +y_i  -z_i+I+{\dfrac{1}{N}}\sum_{j=1}^N  c_{ij}(x_i-x_j),\\
      \dot y_i=c+m(x_i)-y_i, \\
       \dot z_i=\varepsilon(k(x_i-x_0)-z_i).
      \end{array}
\end{equation*}
Setting also  $c=(c_i)$ and  $D=(d_{ij}),$ with 
 $ {\dfrac{1}{N}}\sum\limits_{j=1}^N c_{ij}=c_i,$   $
 d_{ij}=  c_i\delta_{ij}-{\dfrac{1}{N}}c_{ij},
 $  
one can then rewrite the equation for $x_i$ as  
 $$
  \dot x_i= l(x_i) +y_i  -z_i+I+  \sum_{j=1}^N  d_{ij}x_j.
  $$
Defining then
$ {l}({w})=[l(w_1),\dots,l(w_N)]^T ,$
   $ {m}({w})=[m(w_1),\dots,m(w_N)]^T,$
 where functions $l$ and $m$ are defined in \eqref{deflm}, the whole system 
 can be written in vector notation  as
 
 \begin{eqnarray}
 \label{sys4}
 \left [\begin{array}{c}
\dot{{x}}   \\ 
 \dot{{y}}   \\ 
 \dot{{z} }  \\ 
\end{array} \right ]=  \left [\begin{array}{c }
{l}({x}) +{y} -{z} +{I}{e} +{D}{x}     \\ 
{m}( {x})    -{y} +c{e}  \\ 
\varepsilon k{x} -\varepsilon {z} -\varepsilon kx_0{e} \\ 
\end{array} \right ] 
&=&
 \left [\begin{array}{ccc}
 {D} & {\mathbb{I}_N}&- {\mathbb{I}_N}\\
 {\mathbb{O}_N} & - {\mathbb{I}_N} &  {\mathbb{O}_N}\\
 \varepsilon k  {\mathbb{I}_N}&  {\mathbb{O}_N}& -\varepsilon {\mathbb{I}_N}  \\
 \end{array} \right ]
\left[\begin{array}{c}
{x}   \\ 
 {y}   \\ 
 {z}   \\ 
\end{array} \right ] 
\nonumber \\
&+& \left [\begin{array}{c}
{l}({x}) +I{e}  \\ 
 {m}( {x})   + c{e}  \\ 
 -\varepsilon kx_0{e}  \\ 
\end{array} \right ]
 \end{eqnarray} 
 As a consequence, the Jacobian of the vector field on the right hand side is given by
 \begin{equation} 
 \label{j_sys4}
  {J}=    \left [\begin{array}{ccc}
  {\rm diag}({l}^{\prime}({x}))+ {D} & {\mathbb{I}_N}&-{\mathbb{I}_N}\\
 {\rm diag}({m}^{\prime}({x}))& -{\mathbb{I}_N} & \mathbb{O}_N\\
 \varepsilon k {\mathbb{I}_N}& {\mathbb{O}_N}&  -\varepsilon {\mathbb{I}_N}  \\
 \end{array} \right ].
 \end{equation}

  \section{Derivation of efficient implicit solvers}
 \label{sec:solvers} \indent
 The basic idea of our approach will now be presented in the simplest case of the implicit Euler method. 
 The extension to more accurate
 multi-stage or multi-step implicit methods is then straightforward,  see the discussion in Section \ref{sec:esdirk}.  The proposed techniques exploits the special structure of the Jacobians associated to network systems like those described previously. This structure allows to effectively reduce the size of the linear system solved at each iteration of the Newton method when implementing an implicit ODE method. More specifically, at each iteration a linear system is solved by direct methods whose size is equal to that of one
 of the vector unknowns only, rather than to a multiple of it, as it would happen in the case of straightforward application
 of implicit discretizations. The solution of the full system is then recovered by backward substitution in analytically derived expressions, thus reducing the computational cost. The convergence control for the Newton method is carried out on the full solution vector.
 More specifically, the $l^{\infty}$ norm of the increment, normalized by the same norm of the solution at the previous iteration, is taken as an estimate of the relative error. The Newton method employs a tolerance that is much smaller than all the tolerance values used for the time step adaptation techniques (see the discussion in Section \ref{sec:esdirk}).
 All the other components of the solver are identical, so that the proposed simplification does not direct affect the accuracy of the modified methods.

 \subsection{FitzHugh-Nagumo neuron network}
   \label{ssec:sol_FH_nn}
 Consider system \eqref{sys1} and apply the implicit Euler method to compute its numerical solution with time step $h.$
 At each time step, one needs to solve a nonlinear system
 \begin{eqnarray*}
 {G}({u})&=&
  \left [\begin{array}{c }
{G}_1({x}^{n+1},{y}^{n+1})  \\ 
  {G}_2({x}^{n+1},{y}^{n+1})  
\end{array} \right ] 
= \left[\begin{array}{c}
{x}^{n+1}   -h({f}({x}^{n+1}) -{y}^{n+1} + {D}{x}^{n+1}  ) -{x}^n \\ 
 {y}^{n+1} -h( \varepsilon {x}^{n+1}    +  \varepsilon  {g}({y}^{n+1})  )-{y}^n 
 \end{array} \right ]
 = {0},\nonumber
 \end{eqnarray*}
where
$  {u} =\left[ 
{x}^{n+1},
 {y}^{n+1}\right ]^T.
 $
Solving such a system by the Newton method requires to solve linear  systems of the form
 
 $$ {J}_{ G}({u}^{(k)} )\delta^{(k+1)}= \left[{\mathbb{I}_{2N}}-h{J}({u}^{(k)} ) \right]  \delta^{(k+1)}  = -{G}({u}^{(k)} ), $$
 where the vectors ${u}^{(k)}_1,$ ${u}^{(k)}_2 $ denote iterative approximations of 
 ${x}^{n+1}, $ ${y}^{n+1},$ 
 respectively,
 $ \delta^{(k+1)} =\left[ 
 \delta^{(k+1)}_1, \delta^{(k+1)}_2   
   \right ]^T $
 denotes the increment vector such that $ {u}^{(k+1)} = {u}^{(k)} + \delta^{(k+1)} $
 and the Jacobian matrix ${J} $ is the same defined in \eqref{j_sys1}.
Due to the special structure of  $ {J}, $ this system can be rewritten
 as 
 \begin{eqnarray}
   {M}^{(k)}  \delta^{(k+1)}_1    +h  \delta^{(k+1)}_2 
 &=&-{G}_1({u}^{(k)}_1,{u}^{(k)}_2 ),
  \nonumber \\ 
 -h\varepsilon  \delta^{(k+1)}_1 +(1-h\varepsilon a_1) \delta^{(k+1)}_2 &=&
 -{G}_2({u}^{(k)}_1,{u}^{(k)}_2 ),   \nonumber
\end{eqnarray}
where 
 ${M}^{(k)}={\mathbb{I}_N}-h  {D} -h{\rm diag}({f}^{\prime}))$
and the dependency on the quantities at the $k$-th iteration has been omitted for simplicity.
The second equation allows to write, as long as $h\varepsilon a_1<1,$
\begin{equation}
\label{eq:u2_sys0}
 \delta^{(k+1)}_2 =-\frac{1}{1-h\varepsilon a_1}{G}_2 ({u}^{(k)}_1,{u}^{(k)}_2 )
+\frac{h\varepsilon }{1-h\varepsilon a_1}  \delta^{(k+1)}_1.
\end{equation}
As a consequence, the first equation can be rewritten as a linear system whose only unknown is $  \delta{u}^{(k+1)}_1, $ yielding

\begin{eqnarray*}
 \left({M}^{(k)}+\dfrac{\varepsilon h^2}{1-h\varepsilon a_1}{\mathbb{I}_N}\right)  \delta^{(k+1)}_1 
&=&\left( \frac{1-h\varepsilon a_1+\varepsilon h^2}{1-h\varepsilon a_1}{\mathbb{I}_N}-h  {D} -h{\rm diag}({f}^{\prime}) \right)
  \delta^{(k+1)}_1\nonumber \\
&=&-{G}_1({u}^{(k)}_1,{u}^{(k)}_2 )+\dfrac{h}{1-h\varepsilon a_1}{G}_2({u}^{(k)}_1,{u}^{(k)}_2 ).
\end{eqnarray*}
After computing $  \delta{u}^{(k+1)}_1, $   $ \delta{u}^{(k+1)}_2 $ can be obtained substituting back in \eqref{eq:u2_sys0}.

 \subsection{Intracellular Calcium Concentration neuron network}
   \label{ssec:sol_ICC_nn}
 Considering now system \eqref{sys2},
 at each time step one needs to solve  a nonlinear system
 \begin{eqnarray*}
 {G}({u})&=&
  \left [\begin{array}{c }
{G}_1({x}^{n+1},{y}^{n+1},{z}^{n+1})  \\ 
  {G}_2({x}^{n+1},{y}^{n+1},{z}^{n+1}) \\
  {G}_3({x}^{n+1},{y}^{n+1},{z}^{n+1})  
\end{array} \right ] 
= \left[\begin{array}{c}
{x}^{n+1}    \\ 
 {y}^{n+1} \\
 {z}^{n+1}
 \end{array} \right ]-h
  \tau \left [\begin{array}{c }
{f}({x}^{n+1} ) -{y}^{n+1}  -  {\phi}_f({z}^{n+1} )  \\ 
 \varepsilon{K}{D}{x}^{n+1}     +  \varepsilon{K} {g}({y}^{n+1} )  \\ 
\varepsilon {\phi}_r({x}^{n+1} )  +\varepsilon{r}({z}^{n+1} )     \\ 
\end{array} \right ] 
 - \left[\begin{array}{c}
{x}^n    \\ 
 {y}^n \\
 {z}^{n}
 \end{array} \right ]={0},\nonumber
 \end{eqnarray*}
where
$ {u} =\left[ 
{x}^{n+1},
 {y}^{n+1},
 {z}^{n+1}
 \right ]^T.
 $
Solving  
 such a system by the Newton method requires to solve linear  systems of the form
 $$ {J}_{ G}({u}^{(k)} )\delta^{(k+1)}= \left[{\mathbb{I}_{3N}}-h{J}({u}^{(k)} ) \right]  \delta^{(k+1)}  = -{G}({u}^{(k)} ), $$
 where the vectors ${u}^{(k)}_1, $ ${u}^{(k)}_2, $ ${u}^{(k)}_3 $   denote iterative approximations of 
 ${x}^{n+1}, $ ${y}^{n+1},$  ${z}^{n+1},$ 
 respectively,
 $ \delta^{(k+1)} =\left[ 
 \delta^{(k+1)}_1  , 
 \delta^{(k+1)}_2 ,
 \delta^{(k+1)}_3
\right ]^T$  
 denotes the increment vector such that $ {u}^{(k+1)} = {u}^{(k)} +\delta^{(k+1)} $
 and the Jacobian matrix ${J} $ is the same defined in \eqref{j_sys2}.
  Due to the special structure of  ${J}, $ this system can be rewritten
 as
 \begin{eqnarray}
   \left[{\mathbb{I}_N}-\tilde h {\rm diag}({f}^{\prime})\right] \delta^{(k+1)}_1    +\tilde h\delta^{(k+1)}_2 
+\tilde h{\rm diag}({\phi}^{\prime}_f)\delta^{(k+1)}_3 &=&-{G}_1
  \nonumber \\ 
 -\tilde h\varepsilon{K}{D}   \delta^{(k+1)}_1 +({\mathbb{I}_N}-\tilde h\varepsilon a_1{K})
 \delta^{(k+1)}_2 &=&
 -{G}_2  \nonumber\\
 -\tilde h \varepsilon {\rm diag}({\phi}^{\prime}_r) \delta^{(k+1)}_1 +\left( 1+ \frac{\tilde h\varepsilon}{\tau_z}\right)
\delta^{(k+1)}_3 &=& -{G}_3,  \nonumber
\end{eqnarray}
 where we have set $\tilde h=\tau h$ and  the dependency on the quantities at the $k$-th iteration has been omitted for simplicity.
The third equation allows to write
\begin{equation}
\label{eq:u3_sys1}
\delta^{(k+1)}_3 =-\frac{1}{ \left(1+ \frac{\tilde h\varepsilon}{\tau_z}\right)}{G}_3 
+\frac{\tilde h\varepsilon}{\left(1+ \frac{\tilde h\varepsilon}{\tau_z}\right)} {\rm diag}({\phi}^{\prime}_r)\delta^{(k+1)}_1, 
\end{equation}
which can in turn be substituted into the first equation to yield

\begin{eqnarray*}
 \delta^{(k+1)}_2 &=&-\frac{{G}_1 }{\tilde h}+\frac{{\rm diag}( {\phi}^{\prime}_f)}{ \left(1+ \frac{\tilde h\varepsilon}{\tau_z}\right)}{G}_3
 - \left[\frac{{\mathbb{I}_N}}{\tilde h}-{\rm diag}({f}^{\prime})\right] \delta^{(k+1)}_1  \nonumber \\
&-&\frac{\tilde h\varepsilon}{\left(1+ \frac{\tilde h \varepsilon}{\tau_z}\right)} {\rm diag}( {\phi}^{\prime}_f){\rm diag}( {\phi}^{\prime}_r)  \delta^{(k+1)}_1 \nonumber \\
&=&-\frac{{G}_1 }{\tilde h}+\frac{{\rm diag}( {\phi}^{\prime}_f)}{ \left(1+ \frac{\tilde h\varepsilon}{\tau_z}\right)}{G}_3
- {M}^{(k)} \delta^{(k+1)}_1,
 \end{eqnarray*}
 where ${M}^{(k)} $ denotes the diagonal matrix
 $${M}^{(k)}=\left[\frac{{\mathbb{I}_N}}{\tilde h}-{\rm diag}({f}^{\prime})\right]
 +\frac{\tilde h\varepsilon}{\left(1+ \frac{\tilde h \varepsilon}{\tau_z}\right)} {\rm diag}( {\phi}^{\prime}_f){\rm diag}( {\phi}^{\prime}_r).
  $$
As a consequence, the second equation can be rewritten as a linear system whose only unknown is    $  \delta^{(k+1)}_1, $ more specifically
\begin{eqnarray}
\label{eq:u2_sys1}
&&-\tilde h\varepsilon{K}{D}   \delta^{(k+1)}_1 
-\left({\mathbb{I}_N}-\tilde h\varepsilon a_1{K}\right){M}^{(k)} \delta^{(k+1)}_1 \nonumber \\
&&= -{G}_2  +({\mathbb{I}_N}-\tilde h\varepsilon a_1{K})\left[ \frac{{G}_1}{\tilde h} -\frac{{\rm diag}({\phi}^{\prime}_f)}{\left(1+ \frac{\tilde h \varepsilon}{\tau_z}\right)}{G}_3 \right], 
 \end{eqnarray}
 which can be rewritten as
 \begin{equation*}
\left[  ({\mathbb{I}_N}-\tilde h\varepsilon a_1{K}){M}^{(k)}+ \tilde h\varepsilon {K}{D} \right]\delta^{(k+1)}_1 
= {G}_2  +({\mathbb{I}_N}-\tilde h\varepsilon a_1{K})\left[ -\frac{{G}_1 }{\tilde h}
+\frac{{\rm diag}({\phi}^{\prime}_f)}{ \left(1+\frac{\tilde h\varepsilon}{\tau_z}\right)}{G}_3 \right]. 
 \end{equation*}
 After solving this system, the variables $\delta^{(k+1)}_2, \delta ^{(k+1)}_3$ can be calculated substituting back the value
 of  $\delta^{(k+1)}_1 $ in \eqref{eq:u2_sys1}
 and \eqref{eq:u3_sys1}.  
 
   \subsection{Hindmarsh-Rose neuron network}
   \label{ssec:sol_HR_nn}
For  system \eqref{sys4}, one can proceed as in  Section \ref{ssec:sol_ICC_nn}.
   Recalling that in this case $J $
   is given by \eqref{j_sys4},
the nonlinear  systems to be solved at each time step are given by

\begin{eqnarray*}
 &&[{\mathbb{I}_N}- h ({\rm diag}({l}^{\prime})+ {D})] \delta^{(k+1)}_1 - h  \delta^{(k+1)}_2
 + h  \delta^{(k)}_3=-{G}_1\nonumber \\
 &&-h {\rm diag}({m}^{\prime})\delta^{(k+1)}_1  +(1+h)\delta^{(k+1)}_2=-{G}_2\\
&&-h\varepsilon k\delta^{(k+1)}_1 +(1+h\varepsilon)\delta^{(k+1)}_3
  =-{G}_3 \nonumber 
  \end{eqnarray*}
This entails that
$$
 \delta^{(k+1)}_2=-\frac{{G}_2}{(1+h)}+\frac{h}{1+h} {\rm diag}({m}^{\prime})
 \delta^{(k+1)}_1  
\ \ \ 
 \delta^{(k+1)}_3
  =-\frac{{G}_3}{1+h\varepsilon}+\frac{h\varepsilon k}{1+h\varepsilon } \delta^{(k+!)}_1
$$
from which one obtains
\begin{eqnarray*}
&&\left\{[{\mathbb{I}_N}- h ({\rm diag}({l}^{\prime})+ {D})] - \frac{h^2}{1+h} {\rm diag}({m}^{\prime})
+\frac{h^2\varepsilon k}{1+h\varepsilon }{\mathbb{I}_N}\right\} \delta^{(k+1)}_1\nonumber \\
&& =-{G}_1 -\frac{h}{(1+h)}{G}_2 +\frac{h}{1+h\varepsilon}{G}_3.
\end{eqnarray*}

  \section{High order ESDIRK solvers}
 \label{sec:esdirk} \indent

 The procedure
 outlined in  Section \ref{sec:solvers} for the implicit Euler method can be extended in principle to any class of
 higher order implicit methods.
 Here we focus on the specific class of so-called ESDIRK methods (Diagonally Implicit Runge Kutta methods with Explicit first stage, see e.g.
  \cite{hairer:1993, lambert:1991}  for the general terminology on ODE solvers and all the related basic concepts).
As discussed in detail  in \cite{kennedy:2016, kennedy:2019}, to which we refer the reader for the full description of these methods,
the assumption of an explicit first stage   allows to build  in a simpler way stiffly accurate methods which are guaranteed to maintain order two in each intermediate stage. Furthermore, the simplest ESDIRK   is the very widely used TR-BDF2 method, which is in several respects an optimal second order stiff solver \cite{hosea:1996} and is included in the comparisons carried out in this Section. Notice, however, that the choice of ESDIRK methods is mostly due to practical convenience and to concentrate the discussion on a sufficiently broad and relevant class of methods. A reformulation analogous to the one we propose can be carried out for any implicit ODE method and we expect similar advantages to arise for other methods as well.
 We have implemented
 the following methods:
 \begin{itemize}
 \item Second order: ESDIRK2(1)3L[2]SA (section 4.1.1). 
  \item Third order: ESDIRK3(2)4L[2]SA   (section 5.1.1).
   \item Fourth order: ESDIRK4(3)6L[2]SA (section 7.1.1),
 \end{itemize}
 where the acronyms and section numbers are those of reference  \cite{kennedy:2016}.  In particular, in ESDIRKi(j)kL[m]SA i is the order of the method, j is the order of the embedded error estimator, k means the number of stages, m is the stage-order, L stands for L-stable and SA for Stiffly Accurate.
 
  In each stage
 of these methods the same procedure described in Section \ref{sec:solvers} was followed. While this requires a specific
 implementation for each class of neural network systems, as it will be seen in Section \ref{sec:numerical}, this is more than
 compensated by the reduction in computational cost.   
 
  \section{Numerical experiments}
 \label{sec:numerical} \indent
 All the methods discussed in Section \ref{sec:esdirk} have been implemented in MATLAB,
 both in their standard formulation and employing the specific, economical reformulation 
 outlined for each system in Section \ref{sec:solvers}. All the ESDIRK methods considered are endowed with embedded methods
 of lower order,  also described in detail in \cite{kennedy:2016}.  This allows to perform a standard time step adaptation, using as error estimate the difference between
 each method and its corresponding
 embedded method
  (see e.g. again  \cite{hairer:1993}).
  More specifically, at each time level $n,$  two approximations  $u^{n+1}$ and $\hat u^{n+1}$ of the solution at time $t^{n+1}$ are computed by the ESDIRK method  and by the associated embedded method,  with convergence orders $p$ and $\hat p,$ respectively. We set $q = \min (p , \hat p). $ The error estimator is defined as $ e^{n+1}=u^{n+1} - \hat u^{n+1}$. For each component of this vector, 
 one introduces   the quotient
 \begin{equation}
 \label{eq:stepSizeControl_01}
	  \eta^{n+1}_{i}= \frac{| u^{n+1}_{i} - \hat u^{n+1}_{i} | }{ {\tt rtol} | u^{n+1}_{i} | + {\tt atol}}, 
 \end{equation}
where  $\tt atol$ and $\tt rtol $ denote
 absolute and relative error tolerances,   respectively. One then requires that the inequality
 $\eta =  \max_{i=1,\dots,N} \eta_{i}^{n+1} \leq 1 $ 
 holds. 
An optimal choice of the new time step value  is then given by
$$h_{new} = h_{n}  \eta^{-\frac{1}{q+1}}.$$
If $\eta \leq 1 $ is   satisfied,  the solution is advanced with $u_{n+1}$ and the new step size is chosen as $h_{n+1} = h_{new}$.  
  
For all the benchmarks, reference solutions were computed with the $\tt ode15s$ MATLAB function  
using reference tolerances given by $\tt atol_{ref}=10^{-5}atol, $ $ \tt rtol_{ref}=10^{-5}rtol, $ respectively,  where $\tt atol, rtol$ denote the tolerances used by the other solvers.
We have solved the systems  using:
\begin{itemize}
\item
 The $\tt ode23tb$ MATLAB solver, 
because the method it implements is the so called TR-BDF2 method \cite{bank:1985}, \cite{hosea:1996}, which is
 essentially equivalent to the second order  ESDIRK2(1)3L[2]SA  solver.
\item The standard ESDIRK solvers mentioned above, of orders $i=2,3,4, $
denoted as $ESDIRKi,$ respectively.
 \item  The adapted economical  ESDIRK solvers of corresponding orders,
 denoted as $ESDIRKiXY,$ where $XY$ indicates the specific system under consideration.
 \end{itemize}
 In some cases, we have also compared the efficiency of the proposed solvers with that of  reference MATLAB solvers
 such as $\tt ode45$ and $\tt ode15s, $ used with the same value of the tolerance parameters.
 We have computed the errors of each solution with respect  to the reference one as:
$${\cal E}=\dfrac{\max_{k}\|u-u_{ref}\|_{\infty}}{\max_{k}\|u_{ref}\|_{\infty}},$$
where the maximum is over all the computed time steps, as well as
 the CPU time required by each method, denoted by $T_i, T_{iXY}, $ respectively. 
Moreover, we have computed the CPU time ratios  of the standard solvers with respect to economical ones as
$ R_{T,i}= T_{i}/T_{iXY}. $ Values $R_{T,i} > 1 $ denote a superior efficiency of the economical versions of the ESDIRK solvers.

 \subsection{Validation tests}
   \label{ssec:validation}

   The goal of this first set of numerical experiments is to validate our implementation and to assess the sensitivity to the results to
the error tolerance.  
To do this, we consider  a FN network with $N=100$ cells with a sparse lattice connectivity matrix  and the timescale separation $\varepsilon=0.05.$  The initial conditions have been selected as follows. The first component $x$ of each cell is sampled randomly with uniform distribution 
in the interval $[-2,-1]$, once for all computations. The second component has been chosen as   $y=f(x), $ where $f$ is given by \eqref{def:f}, so that the initial conditions are close to the attracting part of the slow manifold  \cite{fenichel:1979}, in order to ensure that the orbit will be  directly approaching the attracting limit cycle of the system.

\begin{figure}
    \centering
    \includegraphics[width=0.35\linewidth,  angle=90]{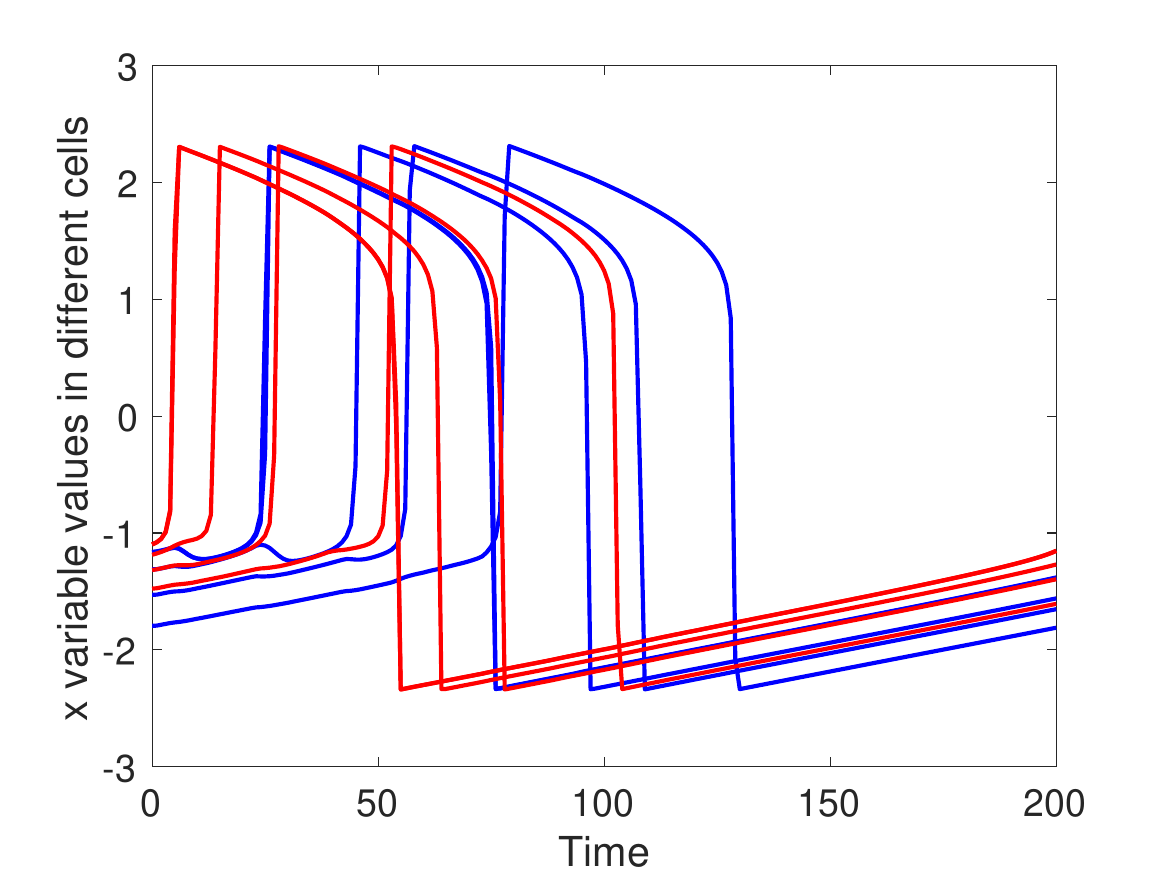} \includegraphics[width=0.35\linewidth, angle=90 ]{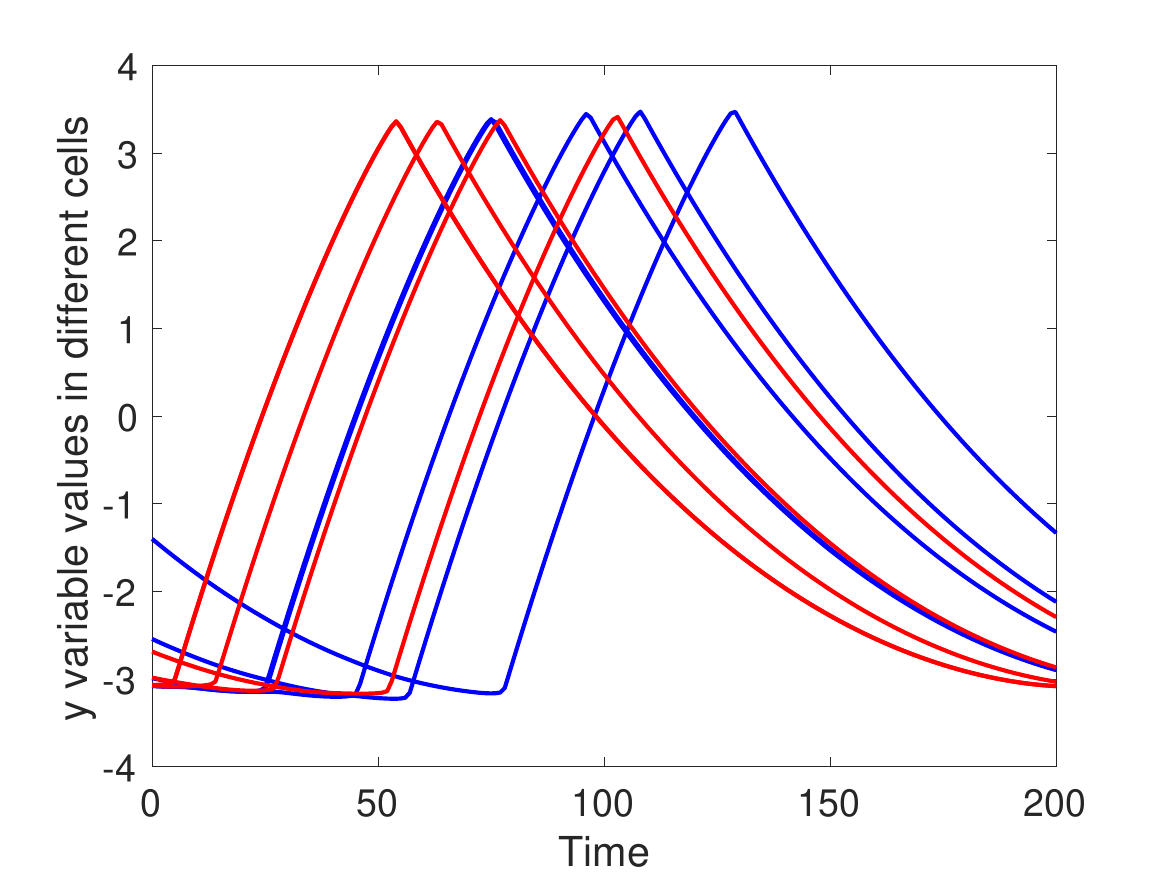}
    \caption{Left panel: reference solution for the $x $ component of the FN system,  cells 1-5 (blue)  and 50-55 (red). Right panel: reference solution for the $y $ component of the FN system,  cells 1-5 (blue)  and 50-55 (red).}
    \label{fig:FN_1}
\end{figure}

In the first test, we have selected a final time $T=200$  and we have varied the absolute and relative tolerances as ${\tt atol}={\tt rtol}=10^{-3},10^{-4}, 10^{-5}.$ Figure \ref{fig:FN_1} shows  the time evolution of the  first component of cells 1-5 and 50-55 for the reference solution in the left panel,
and the time evolution of the  second component of cells 1-5 and 50-55 for the reference solution in the right panel, showing the typical activation and deactivation pattern.

\begin{table}[h!]
    \centering
    \begin{tabular}{ ||c|c|c|c|| }
\hline
   Tolerance   & ESDIRK2   error  & ESDIRK3 
     error  & ESDIRK4 error  \\
    \hline
  $ 10^{-4} $  &  $1.01\times10^{-3}$ &   $1.09\times10^{-3}$&  $5.93\times10^{-4}$ \\ 
    \hline
  $10^{-5}$     &  $8.24\times10^{-5}$  &  $1.76\times10^{-4}$ &  $6.98\times10^{-5}$ \\ 
  \hline
  $10^{-6}$      &  $7.79\times10^{-6}$  &  $1.50\times10^{-5}$ &  $1.90\times10^{-5}$ \\ 
    \hline
    \end{tabular}
    \caption{$l^{\infty}$ relative  errors of the different solvers for varying tolerance values. Errors are computed  for the $x $ component of the first cell in the FN system.}
    \label{tab:errors_test_1}
\end{table}

\begin{table}[h!]
    \centering
    \begin{tabular}{ ||c|c|c|c|| }
\hline
   Tolerance   &     $ R_{T,2}$ &  
    $ R_{T,3}$  &  $ R_{T,4}$ \\
    \hline
  $ 10^{-4} $  & 7.20 & 7.37 & 3.99 \\ 
    \hline
  $10^{-5}$     & 5.61  & 6.31 & 3.94 \\ 
  \hline
  $10^{-6}$     &  4.78  & 5.73  & 4.19\\ 
    \hline
    \end{tabular}
    \caption{CPU time ratios of the different solvers for varying tolerance values in simulation of the FN system.}
    \label{tab:ratios_test_1}
\end{table}

In Table \ref{tab:errors_test_1}, the 
$l^{\infty}$ relative  errors (with respect the reference solution described at the beginning of this Section) of the different economical ESDIRK solvers 
are reported for decreasing tolerance values. Errors are computed  for the $x $ component of the first cell in the FN system, but all component of all cells display similar behaviour.
The errors are of the same order of magnitude of those obtained with the corresponding standard solvers.
In Table \ref{tab:ratios_test_1} we report instead the CPU time ratios  of the standard ESDIRK methods versus their economical counterparts,   highlighting
the superior performance of the proposed reformulation.

\begin{figure}
    \centering
\includegraphics[width=0.65\linewidth]{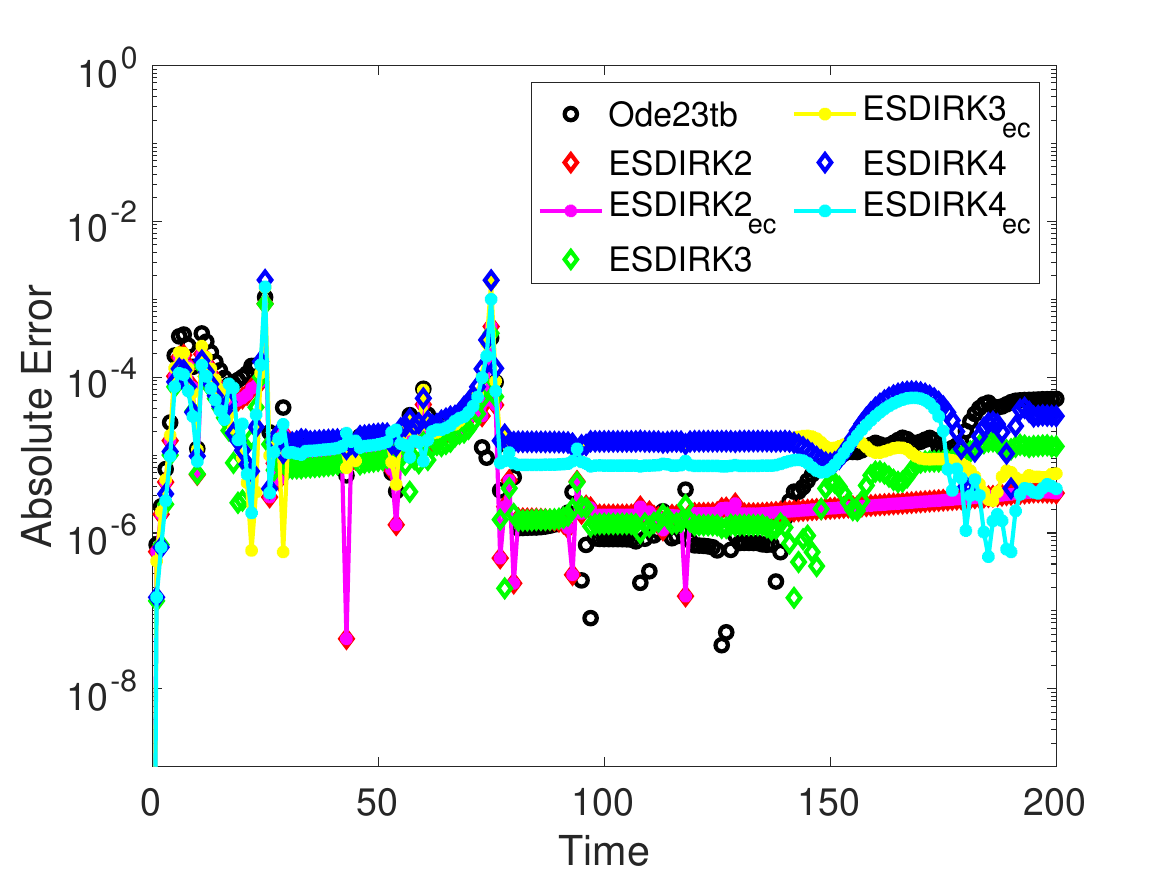}
    \caption{Time evolution of absolute error for $x $ component of the  first cell in the network. Results obtained with tolerance values  ${\tt atol}={\tt rtol}=10^{-4}.$ }
    \label{fig:FN_2}
\end{figure}

A more detailed comparison of the error behaviour of different methods is shown in
 Figure \ref{fig:FN_2}, where we display the absolute  errors for the $x$ component of the  first cell in the network, as computed with tolerance values ${\tt atol}={\tt rtol}=10^{-4}. $ Notice that the sharp peaks in the error evolution correspond to the activation/deactivation phases of the cell kinetics model.

 In a second set of numerical experiments, we study how the performance of the proposed solvers depends on the network size.
   To do this, we consider   again the FN network with the same timescale separation $\varepsilon=0.05$ and the initial conditions selected with the same strategy, assuming $\tt atol=\tt rtol=10^{-4}$ and considering $N=10,20,40, 80, 160, 320 $  cells. The results are shown in Tables \ref{tab:errors_test_2} and \ref{tab:ratios_test_2}. It can be observed that, while the required accuracy is maintained, overall, the economical versions of the solvers are always more efficient than the standard ones, yielding a significant computational cost reduction that is not correlated with the system size.

\begin{table}[h!]
    \centering
    \begin{tabular}{ ||c|c|c|c|| }
\hline
   $N$   & ESDIRK2   error  & ESDIRK3 
     error  & ESDIRK4 error  \\
    \hline
  $ 10 $  &  $1.37\times10^{-3}$ &   $9.48\times10^{-4}$&  $5.52\times10^{-4}$ \\ 
    \hline
  $20$     &  $1.91\times10^{-3}$  &  $2.66\times10^{-3}$ &  $2.74\times10^{-4}$ \\ 
  \hline
  $40$      &  $2.56\times10^{-3}$  &  $1.04\times10^{-2}$ &  $4.01\times10^{-3}$ \\ 
    \hline
    $80$      &  $2.65\times10^{-3}$  &  $2.94\times10^{-3}$ &  $1.38\times10^{-3}$ \\ 
    \hline
    $160$      &  $3.69\times10^{-4}$  &  $8.02\times10^{-4}$ &  $4.38\times10^{-4}$ \\ 
    \hline
    $320$      &  $5.59\times10^{-5}$  &  $5.35\times10^{-5}$ &  $1.76\times10^{-5}$ \\ 
    \hline
    \end{tabular}
    \caption{$l^{\infty}$ relative  errors of the different solvers for varying system size. Errors are computed  for the $x $ component of the first cell in the FN system. The tolerance employed for the adaptive time step choice is $\tt atol=\tt rtol=10^{-4}.$}
    \label{tab:errors_test_2}
\end{table}
\begin{table}[h!]
    \centering
    \begin{tabular}{ ||c|c|c|c|| }
\hline
  $N$   &     $ R_{T,2}$ &  
    $ R_{T,3}$  &  $ R_{T,4}$ \\
    \hline
  $ 10 $  & 5.62 & 2.87 & 2.38 \\ 
    \hline
  $20$     & 7.03  & 4.95 & 2.57 \\ 
  \hline
  $40$     &  7.38  & 6.48  & 3.49\\ 
    \hline
     $80$     &  7.83  & 7.62  & 4.05\\ 
    \hline
     $160$     &  6.60  & 6.45  & 3.68\\ 
    \hline
     $320$     &  5.02  & 4.48  & 2.60\\ 
    \hline
    \end{tabular}
    \caption{CPU time ratios of the different solvers for varying system size in simulation of the FN system.}
    \label{tab:ratios_test_2}
\end{table}
 
 \subsection{Tests with different coupling matrices and timescales}
   \label{ssec:sparse_full}
\begin{figure}
    \centering
    \includegraphics[width=0.30\linewidth]{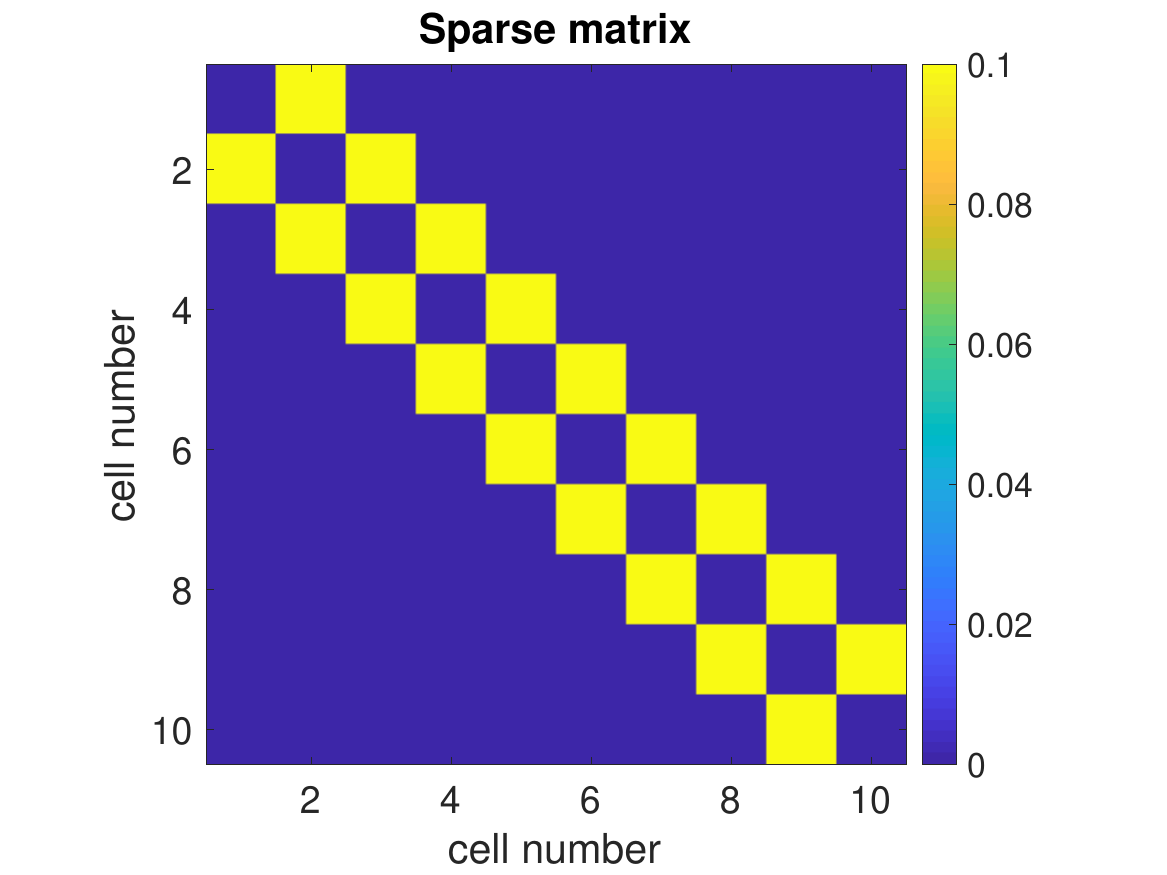}
    \includegraphics[width=0.30\linewidth]{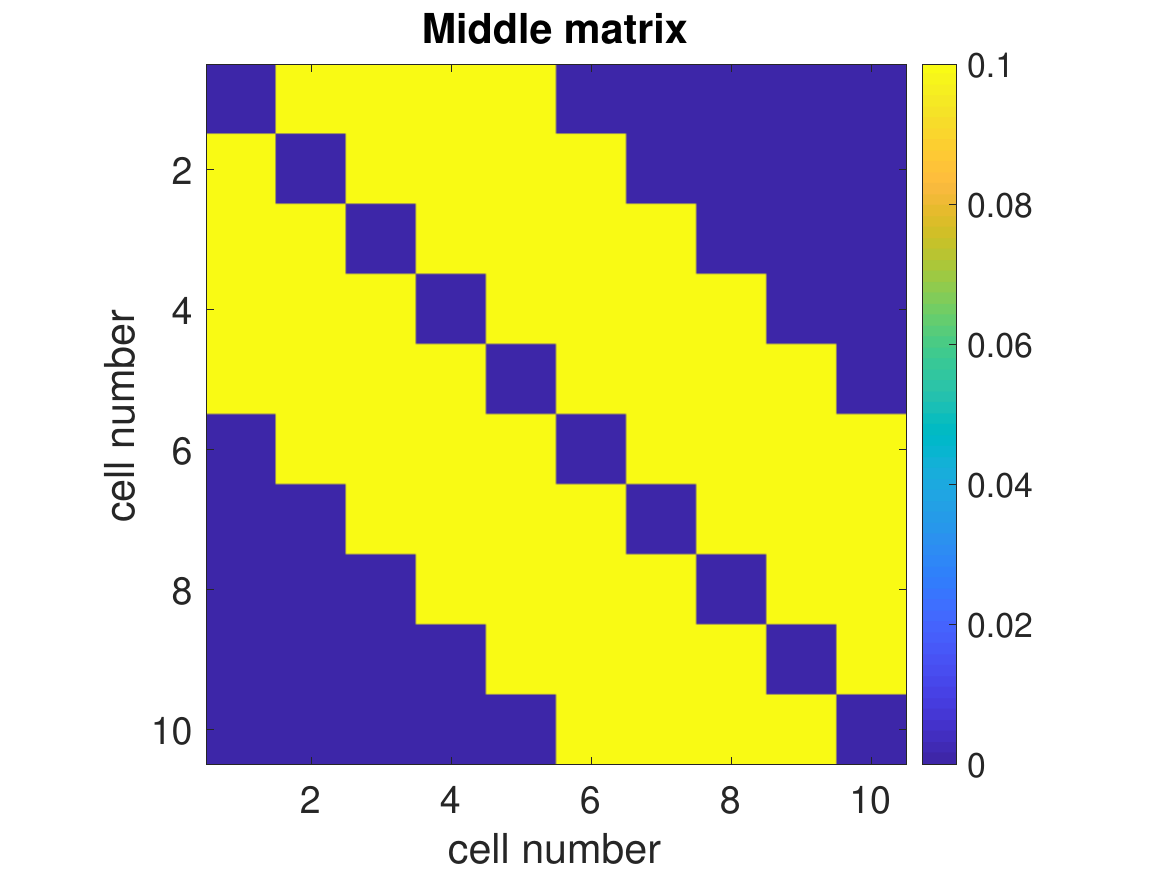}
     \includegraphics[width=0.30\linewidth]{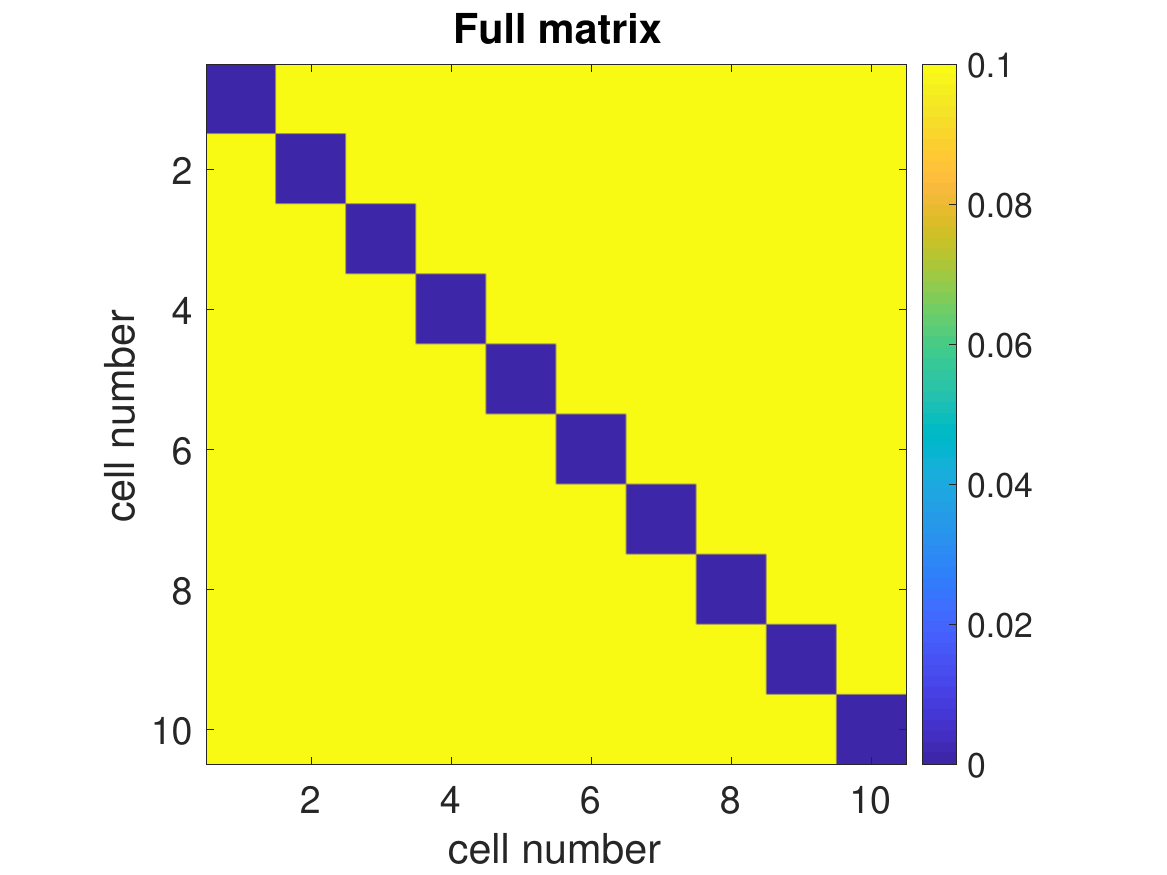}
    \caption{Structure of different coupling matrices for a network of $N=10$ cells.}
    \label{fig_coupling}
 \end{figure}

  In this second set of numerical experiments, we study how the performance depends on the sparsity of the coupling matrix and on the timescale separation.
   To do this, we consider now the HR network with $N=10^3$  cells and  $\tt atol=\tt rtol=10^{-4}$. The initial conditions have been selected in each cell as
$[x_0,y_0,z_0]=[-1.48+\delta,-10.06+\delta,1.84+\delta], $
where $\delta $ denotes a random variable with uniform distribution on $(-0.01,0.01), $   in order to start the simulation in a neighbourhood of the system attractor. 
In the first test, we have selected $\varepsilon=0.01$ and we have considered coupling matrices with the structure represented in Figure \ref{fig_coupling}.  Notice  that, for the sake of clarity, the  plots actually show the case $N=10$ instead of $N=10^3.$  

\begin{figure}
    \centering
    \includegraphics[width=0.55\linewidth,angle=90]{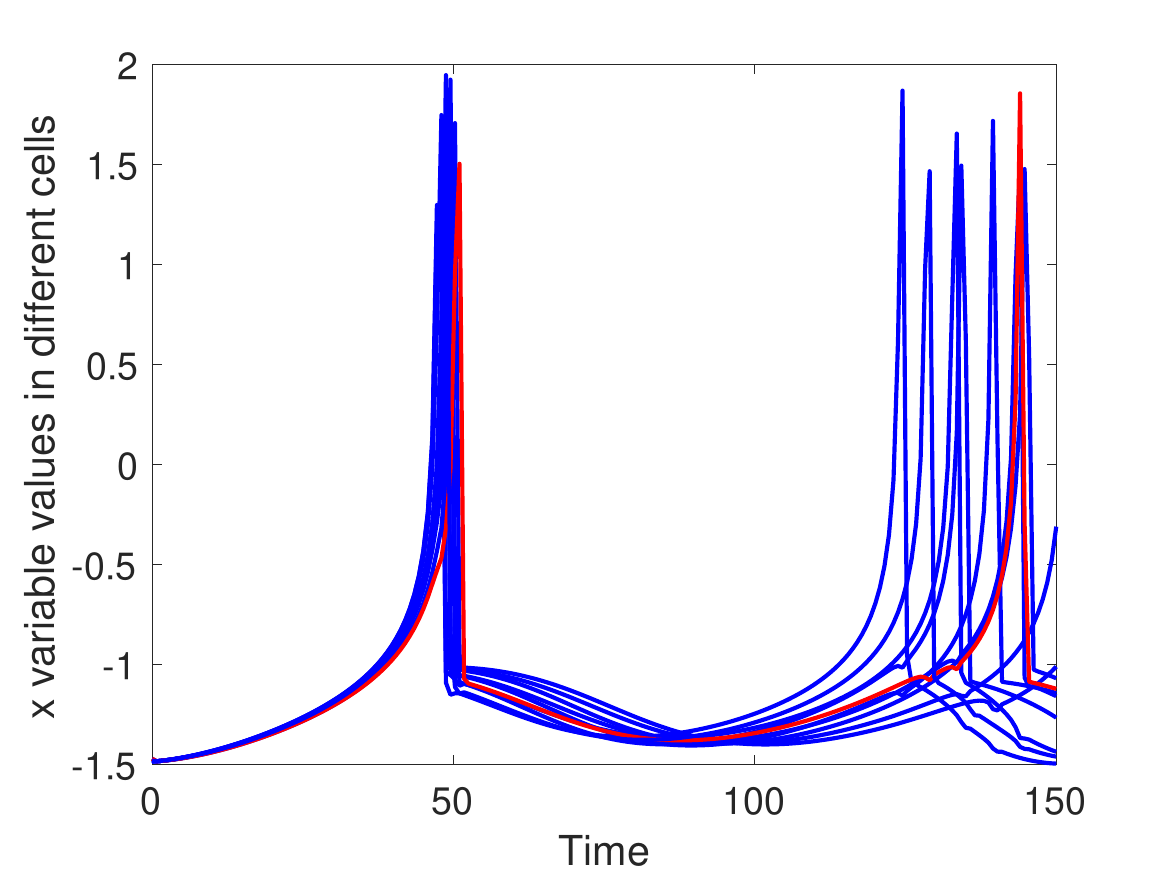}
    \caption{First component of  cells 1-5 (blue) and 501-505 (red) for the reference solution.  }
    \label{fig:HR_sol}
\end{figure}

\begin{figure}
    \centering
    \includegraphics[width=0.65\linewidth]{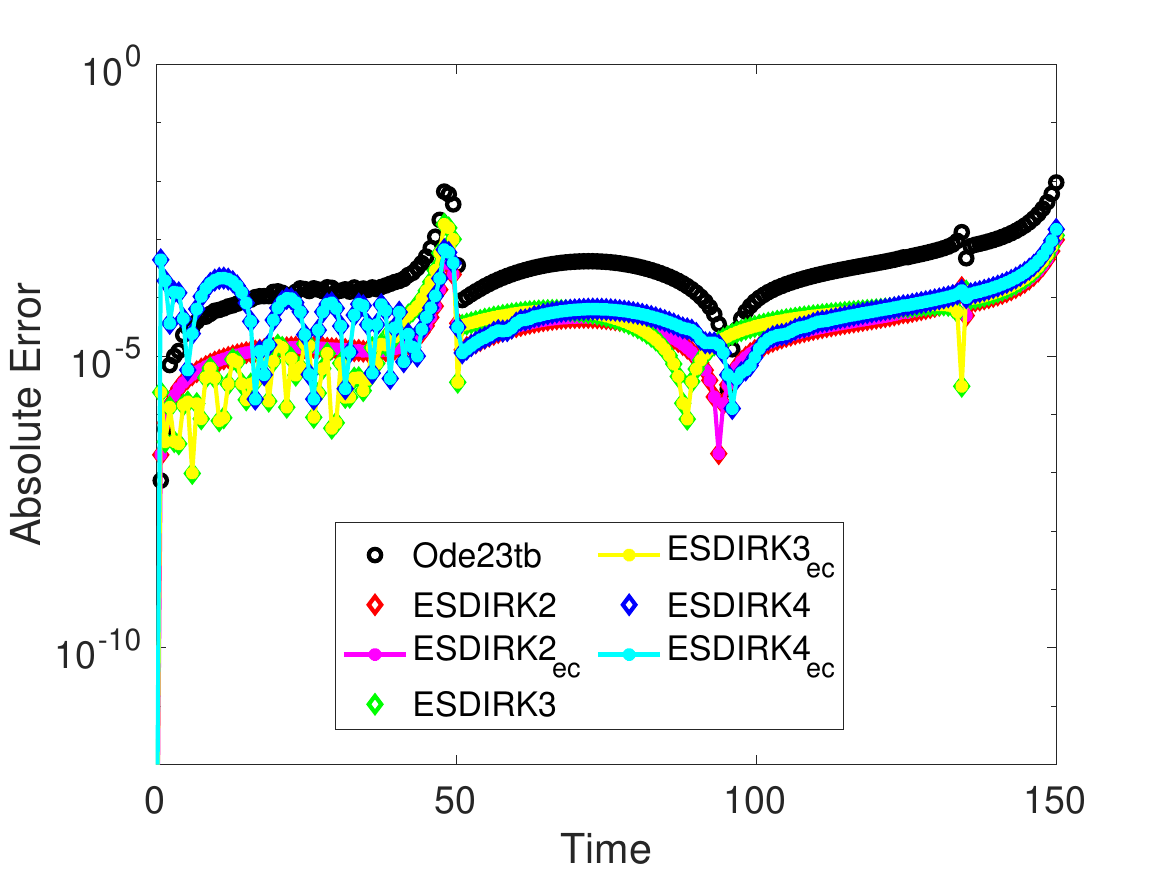} 
        \caption{Time evolution of absolute error for first component, computed with tolerance values  ${\tt atol}={\tt rtol}=10^{-4}.$}
    \label{fig:HR_E}
\end{figure}

Figure \ref{fig:HR_sol} shows the time evolution of the  first component of cells 1-5 and 501-505  for the reference solution in the case of the sparse coupling matrix, displaying an early synchronization peak.
We report the error behaviour of the different methods  in
 Figure \ref{fig:HR_E}, where we display the absolute  errors for the first  component of the  first cell in the network, computed with tolerance values $10^{-4}$. Notice that the sharp peaks in the error evolution correspond to the activation/deactivation phases of the cell kinetics model.

 We then consider  three different types of coupling matrices, in order to assess how they affect the CPU time ratios.
The obtained results are shown in Table \ref{tab:ratios_test_HR1}, where we observe a significant computational cost reduction, specially in the sparse configuration case.

\begin{table}[h!]
    \centering
    \begin{tabular}{ ||c|c|c|c|| }
\hline
   Type of coupling matrix   &     $ R_{T,2}$ &  
    $ R_{T,3}$  &  $ R_{T,4}$ \\
    \hline
  Sparse  & 42.1 & 43.6 & 27.6 \\ 
    \hline
  Middle     & 9.09  & 8.71 & 6.66 \\ 
  \hline
  Full     &  8.13  & 9.86  & 0.05\\ 
    \hline
    \end{tabular}
    \caption{CPU time ratios of the different solvers for different types of coupling matrices in simulation of the HR system.}
    \label{tab:ratios_test_HR1}
\end{table}

 In the second test in this Section, we have  thenconsidered the sparse lattice matrix with $N=10$ cells (see first case in Figure \ref{fig_coupling}), with the same values for the rest of parameters of the fist test in this section, except for the timescale separation parameter, which now has the different values:  $\varepsilon=0.001,0.005,0.01$ and $0.05$. We can see the results in Table \ref{tab:ratios_test_HR2}.   Again, a significant reduction of computational cost
 can be observed, that is not   correlated with the timescale separation variable $\varepsilon$.
  
 \begin{table}[h!]
    \centering
    \begin{tabular}{ ||c|c|c|c|| }
\hline
   $\varepsilon$    &     $ R_{T,2}$ &  
    $ R_{T,3}$  &  $ R_{T,4}$ \\
    \hline
  $ 0.001 $  &  $7.77$ &   $7.25$&  $4.86$ \\ 
    \hline
  $0.005$     &  $8.50$  &  $7.76$ &  $5.33$ \\ 
  \hline
  $0.01$      &  $8.09$  &  $7.85$ &  $5.11$ \\ 
    \hline
    $0.05$      &  $9.45$  &  $3.65$ &  $1.22$ \\ 
    \hline
    \end{tabular}
    \caption{CPU time ratios of the different solvers for different $\varepsilon$ values in simulation of the HR system.}
    \label{tab:ratios_test_HR2}
\end{table}
   
   \subsection{Test with a realistic network configuration}
   \label{ssec:realistic}
\begin{figure}
    \centering
    \includegraphics[width=0.45\linewidth]{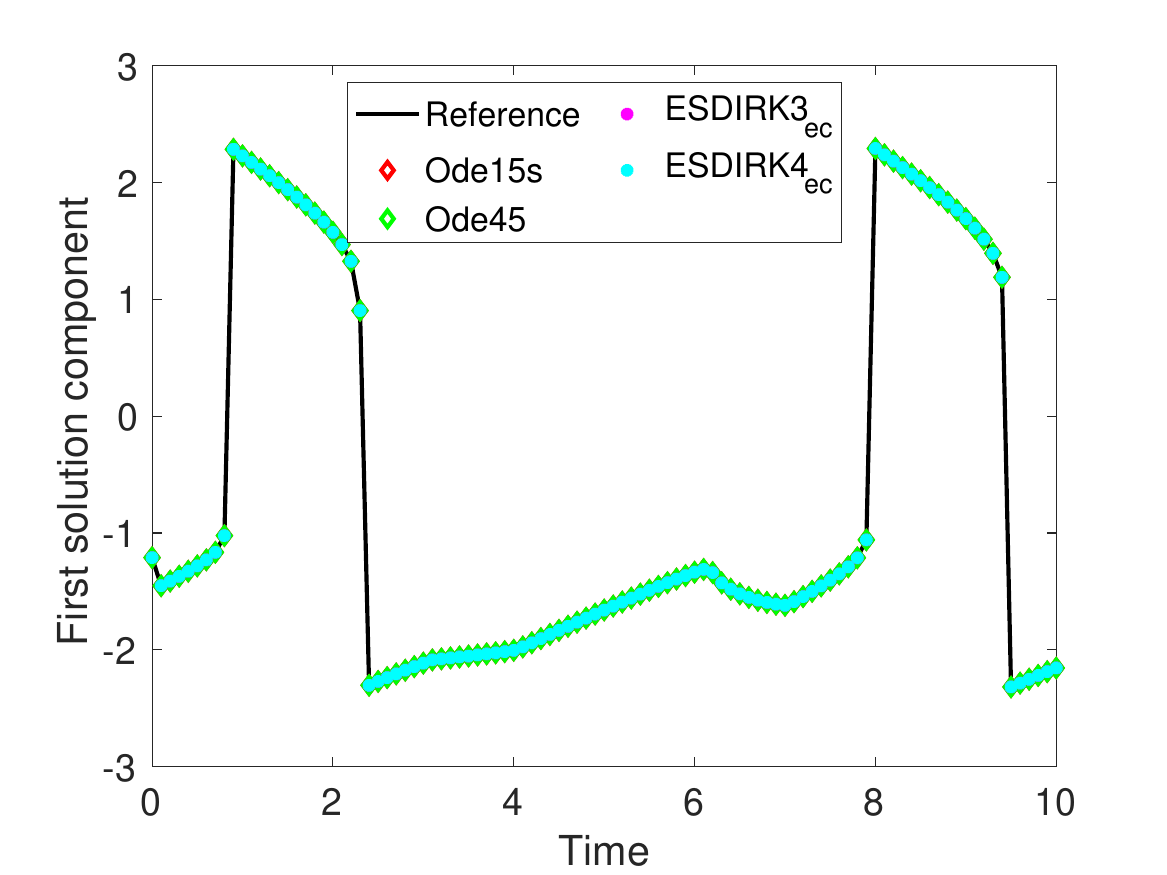}
    \includegraphics[width=0.45\linewidth]{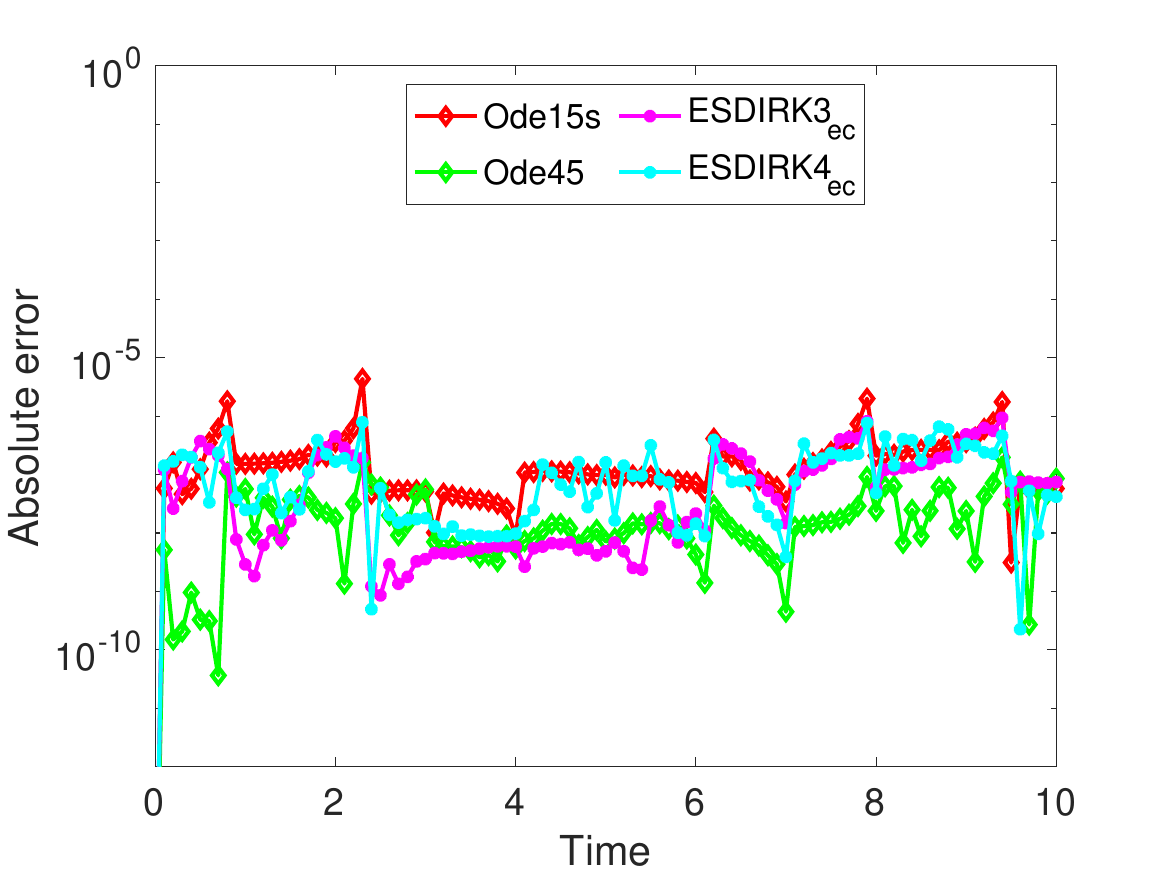}
    \caption{Time evolution of the $x$ variable in the first cell of the first cluster (left plot) and absolute errors for the same variable
    for different solvers (right plot).}
    \label{fig:ICC_1}
 \end{figure}

   We consider here the ICC model on the type of network already studied in \cite{bandera:2022}. 
   The network is   composed of two different clusters  of size 150 cells each, with the cells connected in-phase within each cluster and anti-phase between them. This configuration appears, for instance, in the motoneurons of the embrionic spinal cord of the zebrafish, see  e.g. the discussion in \cite{fallani:2015}.
   The values chosen for the model parameters are the same as in the third example of Section 4.3 of  \cite{bandera:2022}, assuming that the
   first cluster is homogeneous and the second heterogeneous. In this case, only the third and fourth order economical solvers were considered
   and compared to the MATLAB solvers $\tt ode45$ and $\tt ode15s, $ used with the same value of the tolerance parameters
   $\tt atol=\tt rtol=10^{-7}$. In Figure \ref{fig:ICC_1},  we report the time evolution of the $x$ variable in the first cell of the  first cluster   and of the absolute errors on the same variable
    for different solvers.  Notice that the sharp peaks in the error evolution correspond to the activation/deactivation phases of the cell kinetics model.
     It can be observed that, with the same tolerance, the economical ESDIRK solvers achieve  equivalent accuracy
    with respect to the reference MATLAB solvers, even though due to lack of code optimization the required CPU time for the our implementation
    is about one order of magnitude larger than that of the MATLAB {\tt ode45} solver and twice as large as that of the {\tt ode15s} solver.

   \section{Conclusions and future work}
 \label{sec:conclu} \indent
We have outlined a general method to   build efficient, specific versions of standard implicit ODE solvers  tailored for the  
simulation of neural networks. The specific versions of the ODE solvers proposed here   allow to achieve a significant increase in the efficiency
 of network simulations, by reducing the size of the algebraic system being solved at each time step.  

While we have focused here specifically on Explicit first step, Diagonally Implicit Runge Kutta methods (ESDIRK), 
 similar simplifications can be applied to any
implicit ODE solver. In order to demonstrate the capabilities of the proposed methods,
we have considered networks based on three different slow-fast single-cells models, including the classical FitzHugh-Nagumo (FN) model,
the Intracellular Calcium Concentration (ICC) model and the Hindmarsh-Rose (HR) system  model. The numerical
results obtained in a range of simulations of systems with different size and topology demonstrate the potential of the proposed method
to increase substantially the efficiency of numerical simulations of neural networks. 

In future developments,   we plan more extensive applications of the proposed approach to the study of large scale neural networks of biological interest
and to further improve the efficiency by developing self-adjusting multirate extensions of these numerical methods   along the lines of \cite{bonaventura:2020a}. 

\section*{Acknowledgements} 
  {We thank the two anonymous reviewers
for a very careful reading of the original paper and for very detailed comments that have greatly helped to improve the quality of the revised version.}
This work has been supported by Ministerio de Ciencia, Innovación y Universidades into the project  PID2021-123153OB-C21 (Modelos de Orden Reducido H\'ibridos aplicados a flujos incompresibles y redes neuronales cerebrales).
\bibliographystyle{plain}
\bibliography{BFG_2025}
\end{document}